\documentclass[12 pt]{amsart}
\usepackage[utf8]{inputenc}
\usepackage{gensymb}
\usepackage{amsmath}
\usepackage{amssymb}
\usepackage{amsthm}
\usepackage{enumitem}
\usepackage{hyperref}
\usepackage[left = 0.8 in, right = 0.8 in, top = 1.25 in, bottom = 1.25 in]{geometry}
\usepackage{graphicx}
\usepackage{blindtext}

\usepackage{hyperref}
\hypersetup{
     colorlinks=true,
     linkcolor= black,
     filecolor= black,
     citecolor = black,      
     urlcolor= black,
     }
\makeatletter
\newcommand\footnoteref[1]{\protected@xdef\@thefnmark{\ref{#1}}\@footnotemark}
\makeatother

\usepackage{bbm}

\title{Some Results on Zumkeller Numbers}
\author{Sai Teja Somu, Andrzej Kukla, Duc Van Khanh Tran}
\date{}

\theoremstyle{definition}
\newtheorem{theorem}{Theorem}[section]
\newtheorem{lemma}[theorem]{Lemma}
\newtheorem{corollary}[theorem]{Corollary}
\newtheorem{proposition}[theorem]{Proposition}
\newtheorem{conjecture}[theorem]{Conjecture}
\newtheorem{remark}[theorem]{Remark}

\newtheorem{observation}[theorem]{Observation}

\begin{document}


\maketitle

\begin{abstract}\normalfont
    A positive integer $n$ is said to be a Zumkeller number or an integer-perfect number if the set of its positive divisors can be partitioned into two subsets of equal sums. In this paper, we prove several results regarding Zumkeller numbers. For any positive integer $m$, we prove that there are infinitely many positive integers $n$ for which $n+1,\cdots, n+m$ are all Zumkeller numbers. Additionally, we show that every positive integer greater than $94185$ can be expressed as a sum of two Zumkeller numbers and that all sufficiently large integers can be written as a sum of a Zumkeller number and a practical number. We also show that there are infinitely many positive integers that cannot be expressed as a sum of a Zumkeller number and a square or a prime. 
\end{abstract}

\section{Introduction}\label{sec: introduction}

A perfect number, a well-known concept that has been studied since ancient times (see \cite[Ch. 1]{perfect_num_background}), is a natural number that is equal to the sum of its positive proper divisors. In this article,
by natural numbers, we mean positive integers. As a generalization of perfect numbers, a Zumkeller number, or an integer-perfect number, is a natural number $n$ whose set of positive divisors can be partitioned into two subsets of equal sums. Equivalently, a natural number $n$ is Zumkeller if $\frac{\sigma(n)}{2}$ is a sum of distinct positive divisors of $n$, where $\sigma(n)$ denotes the sum of all positive divisors of $n$. Note that a perfect number $n$ is Zumkeller since we can partition its set of positive divisors into one subset containing only $n$ and the other containing all of the positive proper divisors.

The concept of Zumkeller numbers was initially introduced in 1987 by LeVan (see \cite{int_perf}). It was later popularized by Zumkeller, who published this class of numbers on the On-Line Encyclopedia of Integer Sequences (OEIS) in 2003 (see \cite{oeis_zumkeller}). Since then, several authors have studied Zumkeller numbers. In 2008, Clark et al. presented a few results on Zumkeller numbers and proposed two conjectures (see \cite{clark2008}). Later, Peng and Bhaskara Rao provided several results on Zumkeller numbers and their relations with practical numbers (see \cite{Rao_Zumkeller_num}). Practical numbers are natural numbers $n$ such that any natural number less than or equal to $n$ can be written as a sum of distinct positive divisors of $n$, and they were first introduced by Srinivasan in \cite{Srinivasan}. In \cite{Rao_Zumkeller_num}, Peng and Bhaskara Rao also proved one conjecture and partially settled the second conjecture proposed by Clark et al. \cite{clark2008}. In \cite{Mahanta_Zumkeller_num}, Mahanta et al. completely characterized Zumkeller numbers with two distinct prime factors and gave some bounds on Zumkeller numbers with more than two distinct prime factors. In \cite{jokar_differences}, Jokar provided several results on the differences between Zumkeller numbers and proposed a conjecture regarding the existence of strings of consecutive Zumkeller numbers of arbitrary length, which we prove in Section \ref{sec: Zum_arith}.

A few authors have further generalized the notion of Zumkeller numbers and defined and studied $k-$layered numbers, which are natural numbers whose sets of positive divisors can be partitioned into $k$ subsets of equal sums (see \cite{Mahanta_Zumkeller_num}, \cite{jokar_differences}, and \cite{jokar_k-layered}). While Zumkeller and $k-$layered numbers are interesting in their own rights, there have also been many studies applying these concepts into graph theory (see, for example, \cite{Zumkeller_graph_1}, \cite{Zumkeller_graph_2}, \cite{Zumkeller_graph_3}, \cite{Zumkeller_graph_4}, \cite{Zumkeller_graph_5}, \cite{Zumkeller_graph_6}, \cite{Zumkeller_graph_7}, and \cite{Zumkeller_graph_8}).

In this paper, we focus on studying the Zumkeller numbers. More specifically, we prove several results regarding polynomial representations of Zumkeller numbers and sums involving Zumkeller numbers.

In Section \ref{sec: general_res}, we present some general results on polynomial representations of Zumkeller numbers. In particular, we show that for any polynomial of integral coefficients with a positive leading coefficient $P(n)$, there are either no or infinitely many Zumkeller numbers of the form $P(n)$.

In Section \ref{sec: Zum_arith}, we prove the existence of strings of consecutive Zumkeller numbers of arbitrary length. As a corollary, we also show the existence of infinitely many Zumkeller numbers in any arithmetic progression, following the spirit of what has been done for primes and practical numbers by Dirichlet and Margenstern, respectively (see \cite{prime_arithmetic} and \cite[Theorem 10]{practical_arithmetic}).

In Section \ref{sec: additive_rep_zum}, we prove some results regarding additive representations of natural numbers involving Zumkeller numbers. First, we classify natural numbers which can be written as a sum of two Zumkeller numbers. Notably, we classify all even natural numbers that can be written as a sum of two Zumkeller numbers, providing a result analogous to the Goldbach Conjecture for primes (see \cite{goldbach_conj_empirical}) and a result proven by Melfi for practical numbers (see \cite[Theorem 1]{sum_practical}). Then, we prove some results regarding numbers which can be written as a sum of a Zumkeller number and a practical number. Finally, we also provide some results regarding the representations of a natural number as a sum of a Zumkeller number and a square, inspired by \cite[Section 3]{practical_square}, and as a sum of a Zumkeller number and a prime. In particular, we prove that there are infinitely many natural numbers that cannot be written as a sum of a Zumkeller number and a square or as a sum of a Zumkeller number and a prime.

Finally, in Section \ref{sec: fut_research}, we propose a few possible directions for future research on Zumkeller numbers.
\section{General Results on Polynomial Representations of Zumkeller Numbers}\label{sec: general_res}

For the results in this section, we make use of the following result from \cite{Rao_Zumkeller_num}.

\begin{lemma}\label{lem: relatively_prime_multi}
    If $n$ is a Zumkeller number and $w$ is coprime to $n$, then $nw$ is a Zumkeller number.
\end{lemma}
\begin{proof}
    See \cite[Corollary 5]{Rao_Zumkeller_num} for proof.
\end{proof}

The following proposition gives a class of polynomials with integral coefficients $P(n)$ such that $P(n)$ is a Zumkeller number for all non-negative integer $n$.

\begin{proposition}\label{prop: poly_always_Zum}
    Let $c_0, \cdots, c_k$ be non-negative integers. If $a$ is Zumkeller, $a \mid c_i$ for all $1\leq i \leq k$, and $(a,c_0) = 1$, then
    $$P(n) = ac_0 + ac_1n + \cdots +ac_kn^k$$
    is Zumkeller for all non-negative integer $n$.
\end{proposition}
\begin{proof}
    Assume that $a$ is Zumkeller, $a \mid c_i$ for all $1\leq i \leq k$, and $(a,c_0) = 1$. Since $a \mid c_i$ for all $1 \leq i \leq k$ and $(a,c_0)=1$, for all $n \in \mathbb{Z}_{\geq 0}$, 
    $$\left(a, c_0 + c_1n + \cdots + c_kn^k\right) = 1.$$
    By Lemma \ref{lem: relatively_prime_multi}, since $a$ is Zumkeller,
    $$a\left(c_0 + c_1n + \cdots + c_kn^k\right) = ac_0 + ac_1n + \cdots +ac_kn^k$$
    is Zumkeller for all $n\in \mathbb{Z}_{\geq 0}$.
\end{proof}

\begin{remark}\label{rem: 18n+6_18n+12}
    Note that the class of polynomials given in Proposition \ref{prop: poly_always_Zum} is not the only class of polynomials $P(n)$ such that $P(n)$ is Zumkeller for all $n\in\mathbb{Z}_{\geq 0}$. Consider $18n+6$ and $18n+12$, which are of the form $3\cdot 2^k \cdot m$ for some $k\in \mathbb{N}$ and some natural number $m$ coprime to $3\cdot 2^k$. By \cite[Theorem 11]{Rao_Zumkeller_num}, $3\cdot 2^k$ is a Zumkeller number, and by Lemma \ref{lem: relatively_prime_multi}, $3\cdot 2^k\cdot m$ is a Zumkeller number. Hence, $18n+6$ and $18n+12$ are Zumkeller numbers for all non-negative integer $n$. However, $6$ and $12$ are the only two Zumkeller numbers below $18$, and $18n+6$ and $18n+12$ cannot be written as $a(bn+c)$ where $a$ is Zumkeller, $a\mid b$, and $(a,c)=1$.
\end{remark}

The following proposition shows that for any polynomial $P(n)$ with integral coefficients and a positive leading coefficient, there are either no or infinitely many Zumkeller numbers of the form $P(n)$.

\begin{proposition}\label{prop: Zum_poly_infinite}
    Let $P(n)$ be a polynomial with integral coefficients and a positive leading coefficient. If there exists a natural number $\ell$ such that $P(\ell)$ is Zumkeller, then there are infinitely many Zumkeller numbers of the form $P(n)$.
\end{proposition}
\begin{proof}
    Let $P(n)$ be a polynomial of integral coefficients with a positive leading coefficient. Assume $P(\ell) = p_1^{e_1}\cdots p_k^{e_k}$ is Zumkeller, where $\ell, e_1, \cdots, e_k$ are natural numbers, and $p_1, \cdots, p_k$ are primes. Let $m$ be a natural number of the form $m = \ell+p_1\cdots p_k \,P(\ell) \,r$, where $r$ is a natural number.
    
    Then, $m \equiv \ell \, \left(\text{mod } p_i^{e_i+1}\right)$ for all $1\leq i \leq k$, which implies that $P(m) \equiv P(\ell)\, \left(\text{mod } p_i^{e_i+1}\right)$. Thus, $P(m)$ is divisible by $p_i^{e_i}$ but not by $p_i^{e_i+1}$ for all $1\leq i \leq k$. Hence, $P(m)$ is a multiple of $P(\ell)$, and $\frac{P(m)}{P(\ell)}$ is coprime to $P(\ell)$. This means that $P(m)$ is a coprime multiple of $P(\ell)$. 
    
    Since the leading coefficient of $P(n)$ is positive, for sufficiently large $m$, $P(m)$ is positive. Therefore, by Lemma \ref{lem: relatively_prime_multi}, $P(m)$ is Zumkeller for all sufficiently large $m$ of the form $m = \ell+p_1\cdots p_k \, P(\ell)\, r$, and thus there are infinitely many Zumkeller numbers of the form $P(n)$.
\end{proof}

\begin{remark}
Let $n\in \mathbb{Z}_{\geq 0}$. By previous propositions and remark, a polynomial $P(n)$ with integral coefficients and a positive leading coefficient can be Zumkeller for no $n$, infinitely many $n$, or all $n$. For example, $P(n)=n^2$ is never Zumkeller because the sum of divisors of a perfect square is always odd. The polynomial $P(n)=n^3+1$ is Zumkeller for infinitely many $n$ but not all $n$ because $3^3+1 = 28$ is Zumkeller, but $2^3+1 = 9$ is not Zumkeller. Also, as mentioned in Remark \ref{rem: 18n+6_18n+12}, $P(n) = 18n+6$ is Zumkeller for all $n$.
\end{remark}

\section{Strings of Consecutive Zumkeller Numbers and Zumkeller Numbers in Arithmetic Progressions}\label{sec: Zum_arith}
In this section, we prove the existence of strings of consecutive Zumkeller numbers of arbitrary length. Before proving our main result in this section, we need a few lemmas, and the first two lemmas come from \cite{erdos_paper}.

In the lemmas and theorems below, $p_i$ denotes the $i$-th prime, and
$$A_k := p_k\cdots p_{f(k)},$$
where $f(k)$ is the smallest index such that
$$\sigma(p_k\cdots p_{f(k)}) = (1+p_k)\cdots(1+p_{f(k)}) \geq 2 p_k\cdots p_{f(k)}.$$
Such an $f(k)$ always exists, as the product
$$\prod_{p \text{ prime}} \left(1+\frac{1}{p}\right)$$
diverges by \cite[Chapter 7, Theorem 3]{Knopp} and \cite[Theorem 19]{Hardy_Wright}.

\begin{lemma}\label{lem: erdos_lem_1}
    There is an absolute constant $c$ such that every integer $m > cp_k$ is the distinct sum of primes not less than $p_k$.
\end{lemma}
\begin{proof}
    See \cite[Lemma 1]{erdos_paper} for proof.
\end{proof}

\begin{lemma}\label{lem: erdos_lem_2}
    There exists a natural number $k_0$ such that, for every $k > k_0$,
    $$cp_k < m < \sigma(A_k)-cp_k$$
    implies that $m$ is the distinct sum of divisors of $A_k$.
\end{lemma}
\begin{proof}
    See \cite[Lemma 2]{erdos_paper} for proof.
\end{proof}

\begin{lemma}\label{lem: coprime_Zumkeller}
    There exists a natural number $k_1>k_0$ such that $A_k$ is a Zumkeller number for all $k > k_1$.
\end{lemma}
\begin{proof}
    Let $k>k_0$ and
    $$m = \frac{\sigma(A_k)}{2} = (p_k+1) \frac{(p_{k+1}+1)\cdots(p_{f(k)}+1)}{2}.$$
    For sufficiently large $k$,
    $$\frac{(p_{k+1}+1)\cdots(p_{f(k)}+1)}{2} > c,$$
    where $c$ is the absolute constant in Lemma \ref{lem: erdos_lem_1}. Hence,
    $$\frac{\sigma(A_k)}{2} > c(p_k+1) >cp_k,$$
    and
    $$\frac{\sigma(A_k)}{2} = \sigma(A_k) - \frac{\sigma(A_k)}{2} < \sigma(A_k) - cp_k.$$
    By Lemma \ref{lem: erdos_lem_2}, $\frac{\sigma(A_k)}{2}$ can be written as a sum of distinct divisors of $A_k$. Therefore, $A_k$ is a Zumkeller number for all sufficiently large $k$. In other words, there exists a natural number $k_1>k_0$ such that $A_k$ is a Zumkeller number for all $k > k_1$.
\end{proof}

Now we can use the lemmas above to prove our two main theorems of this section. Note that Theorem \ref{thrm: consecutive_arbitrary_length} proves Conjecture 1.30 in \cite{jokar_differences}.

\begin{theorem}\label{thrm: consecutive_arbitrary_length}
    For any natural number $m$, there exist infinitely many natural numbers $n$ such that $n+1, n+2, \cdots, n+m$ are all Zumkeller numbers.
\end{theorem}
\begin{proof}
    Let $j_1,\cdots, j_m$ be natural numbers such that 
    $j_1 > k_1, j_2>f(j_1), j_3> f(j_2), \cdots, j_m > f(j_{m-1})$. By Lemma \ref{lem: coprime_Zumkeller}, $A_{j_1}, A_{j_2}, \cdots, A_{j_m}$ are Zumkeller numbers. Observe that $A_{j_1}^2, A_{j_2}^2, \cdots, A_{j_m}^2$ are pairwise coprime since $j_i > f(j_{i-1})$ for all $2 \leq i \leq m$. Hence, by the Chinese Remainder Theorem, there exist infinitely many natural numbers $n$ satisfying
    \begin{align*}
        n &\equiv -1 + A_{j_1}\, \left(\text{mod } A_{j_1}^2\right),\\
        n &\equiv -2 + A_{j_2}\, \left(\text{mod } A_{j_2}^2\right),\\
        &\:\;\vdots \\
        n &\equiv -m + A_{j_m}\, \left(\text{mod } A_{j_m}^2\right).
    \end{align*}
    For all $1\leq i\leq m$,
    $$n \equiv -i + A_{j_i}\, \left(\text{mod } A_{j_i}^2\right)$$
    implies
    $$n + i \equiv A_{j_i}\, \left(\text{mod } A_{j_i}^2\right),$$
    which implies that $A_{j_i} \mid (n+i)$ and that $\frac{n+i}{A_{j_i}}$ is a natural number coprime to $A_{j_i}$. Thus, for all $1\leq i \leq m$, since $A_{j_i}$ is a Zumkeller number, and $\frac{n+i}{A_{j_i}}$ is coprime to $A_{j_i}$, by Lemma \ref{lem: relatively_prime_multi},
    $$\frac{n+i}{A_{j_i}}\cdot A_{j_i} = n+i$$
    is a Zumkeller number. Therefore, for any natural number $m$, there exist infinitely many natural numbers $n$ such that $n+1, n+2, \cdots, n+m$ are Zumkeller numbers.
\end{proof}
We have the following corollary to Theorem \ref{thrm: consecutive_arbitrary_length}.
\begin{corollary}
Let $a$ be a natural number and $b$ be any non-negative integer. There are infinitely many natural numbers $n$ such that $an+b$ is a Zumkeller number.
\end{corollary}
\begin{proof}
    Let $a$ be a natural number and $b$ be a non-negative integer. For the sake of contradiction, assume that there are only finitely many Zumkeller numbers of the form $an+b$. Then, there exists a natural number $N\geq a+b$ such that no numbers of the form $an+b$ greater than $N$ are Zumkeller numbers. By Theorem \ref{thrm: consecutive_arbitrary_length}, there exists a natural number $k$ such that $k+1, k+2, \cdots, k+2N$ are Zumkeller numbers. At least one of $k+N+1, k+N+2, \cdots, k+2N$ is of the form $an+b$ and is greater than $N$. So, there exists a Zumkeller number of the form $an+b$ greater than $N$, which is a contradiction.
\end{proof}

\section{Additive Representations Involving Zumkeller Numbers}\label{sec: additive_rep_zum}
In this section, we prove some results regarding additive representations of natural numbers involving Zumkeller numbers. First, we classify all natural numbers that are expressible as a sum of two Zumkeller numbers. To do so, we first classify all even natural numbers which can be written as a sum of two Zumkeller numbers and prove that all odd numbers greater than $94185$ are expressible as a sum of two Zumkeller numbers. Then, for odd natural numbers less than or equal to $94185$, we computationally verify whether they can be written as a sum of two Zumkeller numbers.

\begin{lemma}\label{lem: 88m}
    For all natural numbers $m$ such that $11\nmid m$, $88m$ is a Zumkeller number.
\end{lemma}
\begin{proof}
    If $11 \nmid m$, then $88m$ is of the form $2^k\cdot 11 \cdot \ell$, where $k$ and $\ell$ are natural numbers such that $k \geq 3$ and $\left(\ell, 2^k\cdot 11\right) = 1$. By \cite[Theorem 11]{Rao_Zumkeller_num}, $2^k\cdot 11$ is a Zumkeller number for $k\geq 3$. Then, by Lemma \ref{lem: relatively_prime_multi}, since $\ell$ is coprime to $2^k\cdot 11$, $2^k\cdot 11 \cdot \ell$ is a Zumkeller number. Therefore, $88m$ is a Zumkeller number for all natural numbers $m$ such that $11 \nmid m$.
\end{proof}

Now we consider the values of $u_i$ and $v_i$ listed in Table \ref{tab: u_and_v} and make some observations regarding these numbers.

\begin{table}[h]
    \centering
    \bgroup
    \def\arraystretch{1.35}
    \begin{tabular}{lll}
        \hline $i$ & $u_i$ & $v_i$ \\
        \hline 1 & 5985 & 61425 \\
        2 & 17955 & 29835\\
        3 & 2205 & 29925\\
        4 & 14175 & 41895\\
        5 & 26145 & 44625\\
        6 & 6435 & 10395\\
        7 & 22365 & 50085\\
        8 & 6615 & 80535\\
        9 & 18585 & 46305\\
        10 & 2835 & 30555\\
        11 & 14805 & 42525\\
        \hline
    \end{tabular}
    \egroup
    \hspace{12 pt}
    \bgroup
    \def\arraystretch{1.35}
    \begin{tabular}{lll}
        \hline $i$ & $u_i$ & $v_i$ \\
        \hline 12 & 26775 & 82215\\
        13 & 38745 & 94185\\
        14 & 22995 & 50715\\
        15 & 7245 & 34965\\
        16 & 19215 & 74655\\
        17 & 3465 & 7425\\
        18 & 15435 & 33915\\
        19 & 8925 & 27405\\
        20 & 11655 & 39375\\
        21 & 23625 & 69825\\
        22 & 7875 & 63315\\
        \hline
    \end{tabular}
    \egroup
    \hspace{12 pt}
    \bgroup
    \def\arraystretch{1.35}
    \begin{tabular}{lll}
         \hline $i$ & $u_i$ & $v_i$ \\
         \hline 23 & 19845 & 84525\\
        24 & 4095 & 31815\\
        25 & 6825 & 16065\\
        26 & 9555 & 28035\\
        27 & 12285 & 67725\\
        28 & 5775 & 8415\\
        29 & 8505 & 36225\\
        30 & 20475 & 48195\\
        31 & 4725 & 23205\\
        32 & 16695 & 25935\\
        33 & 945 & 28665\\
        \hline
    \end{tabular}
    \egroup
    \hspace{12 pt}
    \bgroup
    \def\arraystretch{1.35}
    \begin{tabular}{lll}
        \hline $i$ & $u_i$ & $v_i$ \\
        \hline 34 & 12915 & 31395\\
        35 & 24885 & 34125\\
        36 & 9135 & 36855\\
        37 & 21105 & 39585\\
        38 & 5355 & 33075\\
        39 & 8085 & 17325\\
        40 & 1575 & 29295\\
        41 & 13545 & 33345\\
        42 & 25515 & 53235\\
        43 & 9765 & 37485\\
        44 & 21735 & 58695\\
        \hline
    \end{tabular}
    \egroup
    \vspace{8 pt}
    \caption{The values of $u_i$ and $v_i$}
    \label{tab: u_and_v}
\end{table}

\begin{observation}\label{obs: 2}
    For all $1\leq i \leq 44$, $u_i \equiv v_i\, (\text{mod } 88)$.
\end{observation}

\begin{observation}\label{obs: 3}
    For all $1\leq i \leq 44$, $u_i$ and $v_i$ are Zumkeller numbers.
\end{observation}

\begin{observation}\label{obs: 4}
    Let $\displaystyle{d_i = \frac{v_i-u_i}{88}}$, then $11 \nmid d_i$ for all $1\leq i \leq 44$.
\end{observation}

\begin{lemma}\label{lem: u_and_v}
    All odd numbers $n$ greater than $94185$ can be expressed as
    $$n = 88m + u_i = 88(m-d_i) + v_i$$
    for some $1\leq i \leq 44$, where $m$ is a natural number such that $m > d_i$.
\end{lemma}
\begin{proof}
    Let $n$ be an odd number greater than $94185$. First, note that for all $1\leq i \leq 44$,
    $$u_i \equiv (2i-1)\text{ (mod 88)}.$$
    Since $n$ is odd, $n\equiv u_i\text{ (mod 88)}$ for some $1\leq i \leq 44$. This implies that
    $$n = 88m + u_i = 88(m-d_i) + v_i$$
    for some natural number $m$ and some $1\leq i \leq 44$. For the sake of contradiction, assume $m \leq d_i$. Then,
    $$n = 88(m-d_i) + v_i \leq v_i,$$
    which implies that $n \leq 94185$ because $\max\{u_i, v_i\} = 94185$. This contradicts the assumption that $n > 94185$. Therefore,
    $$n = 88m + u_i = 88(m-d_i) + v_i$$
    for some $1\leq i \leq 44$, where $m > d_i$.
\end{proof}

\begin{proposition}
    Let $n$ be a natural number.
    \begin{enumerate}
        \item[(a)] If $n$ is even, then $n$ can be written as a sum of two Zumkeller numbers if and only if $n \geq 12$ and $n\not\in \{14, 16, 20, 22, 28, 38\}$.
        \item[(b)] If $n$ is an odd number greater than $94185$, then $n$ can be written as a sum of two Zumkeller numbers.
    \end{enumerate}
\end{proposition}
\begin{proof}{\ \\}
    (a) Let $n$ be an even natural number. Then, $n$ is congruent to 0, 2, 4, 6, 8, 10, 12, 14, or 16 modulo 18. From the OEIS entry \cite{oeis_zumkeller}, we know that 6, 12, 20, 28, 40, and 80 are Zumkeller numbers. Recall from Remark \ref{rem: 18n+6_18n+12} that $18m+6$ and $18m+12$ are Zumkeller for all non-negative integers $m$. For $k\geq 5$, we have
    \begin{align*}
        18k &= 18(k-1) + 6 + 12, \\
        18k+2 &= 18(k-5) + 12 + 80,\\
        18k+4 &= 18(k-2) + 12 + 28,\\
        18k+6 &= 18(k-1) + 12 + 12, \\
        18k+8 &= 18(k-1)+6 + 20,\\
        18k+10 &= 18(k-2)+6  +  40, \\
        18k+12 &=  18k+6 + 6,\\
        18k+14 &= 18(k-1)+12 + 20,\\
        18k+16 &= 18(k-2)+12  + 40.
    \end{align*}
    So, all even numbers $n \geq 18\cdot 5 = 90$ can be written as a sum of two Zumkeller numbers. By referring to the OEIS entry \cite{oeis_zumkeller}, we can check that all even numbers $12 \leq n \leq 88$ such that $n\not\in \{14, 16, 20, 22, 28, 38\}$ are expressible as a sum of two Zumkeller numbers:
    $$\begin{tabular}{ccc}
    \:{$\!\begin{aligned}
        12  &=  6   +  6,\\
         
        26  &=  20  +  6,\\
        
        34  &=  28  +  6,\\
        
        42  &=  30  +  12,\\
        
        48  &=  42  +  6,\\
        
        54  &=  48  +  6,\\
        
        60  &=  54  +  6,\\
        
        66  &=  60  +  6,\\
        
        72  &=  66  +  6,\\
        
        78  &=  66  +  12,\\
        
        84  &=  78  +  6,\\
        
    \end{aligned}$}\:    &
    \:{$\!\begin{aligned}
        
        18  &=  12  +  6,\\

        30  &=  24  +  6,\\

        36  &=  30  +  6,\\

        44  &=  24  +  20,\\

        50  &=  30  +  20,\\

        56  &=  28  +  28,\\

        62  &=  56  +  6,\\

        68  &=  56  +  12,\\

        74  &=  54  +  20,\\

        80  &=  60  +  20,\\

        86  &=  80  +  6,\\
        
    \end{aligned}$}\: &
    \:{$\!\begin{aligned}
        
        24  &=  12  +  12,\\ 
        
        32  &=  20  +  12,\\
        
        40  &=  28  +  12,\\
        
        46  &=  40  +  6,\\
        
        52  &=  40  +  12,\\
        
        58  &=  30  +  28,\\
        
        64  &=  40  +  24,\\
        
        70  &=  42  +  28,\\
        
        76  &=  70  +  6,\\
        
        82  &=  70  +  12,\\
        
        88  &=  60  +  28.
    \end{aligned}$}\:
    \end{tabular}$$
    Since the smallest Zumkeller number is 6, all even natural numbers $n < 6+6=12$ are not representable as a sum of two Zumkeller numbers. By referring to the OEIS entry \cite{oeis_zumkeller}, we can also check that $n$ cannot be written as a sum of two Zumkeller numbers if $n\in \{14, 16, 20, 22, 28, 38\}$.\\\\
    (b) Let $n$ be an odd number greater than $94185$. By Lemma \ref{lem: u_and_v},
    $$n = 88m + u_i = 88(m-d_i) + v_i$$
    for some $1\leq i \leq 44$, where $m > d_i$. By Observation \ref{obs: 4},  $11\nmid d_i$ for all $1\leq i \leq 44$, so
    $$m \not\equiv (m - d_i) \, (\text{mod } 11).$$
    This implies that $11 \nmid m$ or $11 \nmid (m-d_i)$. Then, by Lemma \ref{lem: 88m}, $88m$ or $88(m-d_i)$ is a Zumkeller number. By Observation \ref{obs: 3}, $u_i$ and $v_i$ are Zumkeller numbers for all $1\leq i \leq 44$. Therefore, $n$ can be written as a sum of two Zumkeller numbers.
\end{proof}

Since the smallest Zumkeller number is $6$ and the smallest odd Zumkeller number is $945$, any odd natural number less than $945+6=951$ cannot be written as a sum of two Zumkeller numbers. Now we describe the algorithm to verify whether an odd number larger than $951$ and less than or equal to $94185$ is expressible as a sum of two Zumkeller numbers.

In our algorithm, we first generate all Zumkeller numbers below $94185$ and then check whether each natural number less than or equal to $94185$ is expressible as a sum of two Zumkeller numbers. To generate Zumkeller numbers, we make use of the following result.
\begin{proposition}
    A natural number $n$ is Zumkeller if and only if $\frac{\sigma(n)}{2} - n$ can be written as a sum of distinct proper positive divisors of $n$.
\end{proposition}
\begin{proof}
    See \cite[Proposition 3]{Rao_Zumkeller_num} for proof.
\end{proof}
To optimize the generating algorithm, we also make use of the fact that all natural numbers satisfying any of the congruences below are Zumkeller numbers:
$$\begin{tabular}{cccc}
   \:{$\!\begin{aligned}
       &6\, (\text{mod } 18),\\
       &60 \, (\text{mod } 100),\\
       &84 \, (\text{mod } 196),\\
       &88 \, (\text{mod } 968),\\
       &440 \, (\text{mod } 968),\\
       &792 \, (\text{mod } 968),\\
       &312 \, (\text{mod } 1352),\\
       &728 \, (\text{mod } 1352),\\
       &1144 \, (\text{mod } 1352),
   \end{aligned}$}\: &
   \:{$\!\begin{aligned}
       &12 \, (\text{mod } 18),\\
       &80 \, (\text{mod } 100),\\
       &112 \, (\text{mod } 196),\\
       &176 \, (\text{mod } 968),\\
       &528 \, (\text{mod } 968),\\
       &880 \, (\text{mod } 968),\\
       &416 \, (\text{mod } 1352),\\
       &832 \, (\text{mod } 1352),\\
       &1248 \, (\text{mod } 1352).
   \end{aligned}$}\: &
   \:{$\!\begin{aligned}
       &20 \, (\text{mod } 100),\\
       &28 \, (\text{mod } 196),\\
       &140 \, (\text{mod } 196),\\
       &264 \, (\text{mod } 968),\\
       &616 \, (\text{mod } 968),\\
       &104 \, (\text{mod } 1352),\\
       &520 \, (\text{mod } 1352),\\
       &936 \, (\text{mod } 1352),\\
       &
   \end{aligned}$}\: &
   \:{$\!\begin{aligned}
       &40 \, (\text{mod } 100),\\
       &56 \, (\text{mod } 196),\\
       &168 \, (\text{mod } 196),\\
       &352 \, (\text{mod } 968),\\
       &704 \, (\text{mod } 968),\\
       &208 \, (\text{mod } 1352),\\
       &624 \, (\text{mod } 1352),\\
       &1040 \, (\text{mod } 1352),\\
       &
   \end{aligned}$}\:
\end{tabular}$$
All natural numbers satisfying the congruences above can be shown to be Zumkeller by using a similar reasoning to that in Remark \ref{rem: 18n+6_18n+12}, where we showed that all natural numbers congruent to $6$ or $12$ modulo $18$ are Zumkeller numbers.

The algorithm described above is implemented with Python in \verb|NonZumkellerSums.ipynb| on \href{https://github.com/somuteja/ZumkellerNumbers}{GitHub}\footnote{\label{foot: 1}\href{https://github.com/somuteja/ZumkellerNumbers}{https://github.com/somuteja/ZumkellerNumbers}}. As a result, we find that the largest number that cannot be written as a sum of two Zumkeller numbers is $32761$ and that the number of odd numbers between $953$ and $32761$ that cannot be written as a sum of two Zumkeller numbers is $1055$. All natural numbers not expressible as a sum of two Zumkeller numbers are listed in \verb|NonZumkellerSumsList.csv| on \href{https://github.com/somuteja/ZumkellerNumbers}{GitHub}\footnoteref{foot: 1}.

Next, we prove some results regarding numbers which are expressible as a sum of a Zumkeller number and a practical number. Practical numbers, first introduced by Srinivasan in \cite{Srinivasan}, are natural numbers $n$ such that any natural number less than or equal to $n$ can be written as a sum of distinct positive divisors of $n$.

We first give some classes of natural numbers which can be written as a sum of a Zumkeller number and a practical number. Notably, we show that all even numbers greater than or equal to $8$ are expressible as a sum of a Zumkeller number and a practical number.

\begin{proposition}
    Let $n$ be an even natural number or a natural number congruent to $7$ or $13$ modulo $18$. Then, $n$ can be written as a sum of a Zumkeller number and a practical number if and only if $n \geq 7$.
\end{proposition}
\begin{proof}
    Let $n$ be an even natural number or a natural number congruent to $7$ or $13$ modulo $18$. Then, $n$ is congruent to $0$, $2$, $4$, $6$, $7$, $8$, $10$, $12$, $13$, $14$, or $16$ modulo $18$. Recall from Remark \ref{rem: 18n+6_18n+12} that $18m+6$ and $18m+12$ are Zumkeller for all non-negative integers $m$. From the OEIS entry \cite{oeis_practical}, we know that $1$, $2$, $4$, $6$, $8$, $12$, and $16$ are practical numbers. For $k\geq 1$, we have
    \begin{align*}
        18k &= 18(k-1)+ 6 + 12, \\
        18k+2 &= 18(k-1)+12+8,\\
        18k+4 &= 18(k-1)+6+16,\\
        18k+6 &= 18(k-1)+12+12.
    \end{align*}
    For $k\geq 0$, we have
    \begin{align*}
        18k+7 &= 18k+6+1,\\
        18k+8 &=  18k+6+2, \\
        18k+10 &= 18k+6+4, \\
        18k+12 &= 18+6+6,\\
        18k+13 &= 18k+12+1,\\
        18k+14 &= 18k+12+2,\\
        18k+16 &= 18k+12+4.
    \end{align*}
    So, even natural numbers $n$ and natural numbers $n$ congruent to $7$ or $13$ modulo $18$ can be written as a sum of a Zumkeller number and a practical number if $n \geq 7$. Since the smallest Zumkeller number is $6$ and the smallest practical number is $1$, $n$ cannot be written as a sum of a Zumkeller number and a practical number if $n<6+1=7$.
\end{proof}

Now we prove that all sufficiently large integers can be written as a sum of a Zumkeller number and a practical number. In particular, this shows that all sufficiently large odd numbers can be written as a sum of a Zumkeller number and a practical number.

\begin{lemma}\label{lem: practical_arith}
    Let $a, b \in \mathbb{N}$. If $2 \nmid a$, then there are infinitely many $n \in \mathbb{Z}_{\geq 0}$ such that $an + b$ is practical.
\end{lemma}
\begin{proof}
    See \cite[Lemma 2.3]{special_form_practical} for proof.
\end{proof}

\begin{proposition}
    There exists a natural number $N$ such that all natural numbers $n\geq N$ can be written as a sum of a Zumkeller number and a practical number. 
\end{proposition}
\begin{proof}
    Let $n$ be a natural number. Let $i$ be a natural number such that $n-945\equiv i \left(\text{mod } 945^2\right)$ and $1\leq i \leq 945^2$. By Lemma \ref{lem: practical_arith}, for all $1 \leq i \leq 945^2$, there exists a practical number $a_i$ such that
    $a_i \equiv i \, \left(
    \text{mod } 945^2\right).$
    Let $N = \max\{a_i\}_{i=1}^{945^2}$, and suppose $n\geq N$. Then,
    $$n - a_i \geq n-N \geq 0,$$
    and
    $$n - a_i \equiv 945 \left(\text{mod } 945^2\right).$$
    So,
    $$n-a_i = 945^2 k + 945$$
    for some non-negative integer $k$. Since $945$ is Zumkeller, by Proposition \ref{prop: poly_always_Zum}, $945^2 k + 945$ is Zumkeller for all non-negative integer $k$. Therefore, $n-a_i$ is Zumkeller, and thus $n$ can be written as a sum of a Zumkeller number and a practical number.
\end{proof}

Lastly, we provide some results regarding the sums of a Zumkeller number and a square or a prime. First, we give some classes of natural numbers which can be written as a sum of a Zumkeller number and a square.

\begin{proposition}\label{prop: Zum_plus_square_mod_18}
    Let $n$ be a natural number congruent to $1$, $3$, $4$, $6$, $7$, $10$, $12$, $13$, $15$, or $16$ modulo $18$. Then, $n$ can be written as a sum of a Zumkeller number and a square if and only if $n\not\in \{1,3,4,6,12,19,30\}$.
\end{proposition}
\begin{proof}
    Let $n$ be a natural number congruent to $1$, $3$, $4$, $6$, $7$, $10$, $12$, $13$, $15$, or $16$ modulo $18$. Recall from Remark \ref{rem: 18n+6_18n+12} that $18m+6$ and $18m+12$ are Zumkeller for all non-negative integers $m$. For $k\geq 2$, we have
    \begin{align*}
        18k + 1 &= 18(k-2)+12 + 25, \\
        18k+3 &= 18(k-1)+12+9,\\
        18k+4 &= 18(k-1)+6 + 16,\\
        18k+6 &=  18(k-2)+6 + 36, \\
        18k+7 &= 18k+6 + 1,\\
        18k+10 &= 18k+ 6 + 4, \\
        18k+12 &=  18(k-2)+ 12 + 36,\\
        18k+13 &= 18k+12 + 1,\\
        18k+15 &= 18k+6+9,\\
        18k+16 &= 18k+12  + 4.
    \end{align*}
    So, natural numbers $n$ congruent to $1$, $3$, $4$, $6$, $7$, $10$, $12$, $13$, $15$, or $16$ modulo $18$ can be written as a sum of a Zumkeller number and a square if $n \geq 18\cdot 2 + 1 = 37$. By referring to the OEIS entry \cite{oeis_zumkeller}, we can check that all natural numbers $n$ congruent to $1$, $3$, $4$, $6$, $7$, $10$, $12$, $13$, $15$, or $16$ modulo $18$ such that $n<37$ and $n\not\in \{1,3,4,6,12,19,30\}$ are expressible as a sum of a Zumkeller number and a square:
    $$\begin{tabular}{ccccc}
      {$\!\begin{aligned}
          7 &= 6+1,\\
          21 &= 20+1,\\
          31 &= 30+1,
      \end{aligned}$} &
      {$\!\begin{aligned}
          10 &= 6+4,\\
          22 &= 6 + 16,\\
          33 &= 24 + 9, 
      \end{aligned}$} &
      {$\!\begin{aligned}
          13 &= 12+1,\\
          24 &= 20 + 4,\\
          34 &= 30+4.
      \end{aligned}$} &
      {$\!\begin{aligned}
          15 &= 6+9,\\
          25 &= 24 + 1,\\
          &
      \end{aligned}$} &
      {$\!\begin{aligned}
          16 &= 12+4,\\
          28 &= 24 + 4,\\
          &
      \end{aligned}$}
    \end{tabular}$$
    By referring to the OEIS entry \cite{oeis_zumkeller}, we can also check that $n$ cannot be written as a sum of a Zumkeller number and a square if $n\in \{1,3,4,6,12,19,30\}$.
\end{proof}

Now we give some classes of natural numbers expressible as a sum of a Zumkeller number and a prime.

\begin{proposition}\label{prop: Zum_plus_prime_mod_18}
    Let $n$ be a natural number congruent to $1$, $5$, $7$, $8$, $9$, $11$, $13$, $14$, $15$, or $17$ modulo $18$. Then, $n$ can be written as a sum of a Zumkeller number and a prime if and only if $n \geq 8$.
\end{proposition}
\begin{proof}
    Let $n$ be a natural number congruent to $1$, $5$, $7$, $8$, $9$, $11$, $13$, $14$, $15$, or $17$ modulo $18$. Recall from Remark \ref{rem: 18n+6_18n+12} that $18m+6$ and $18m+12$ are Zumkeller for all non-negative integers $m$. For $k\geq 1$, we have
    \begin{align*}
        18k + 1 &= 18(k-1)+6 + 13, \\
        18k+5 &= 18(k-1)+6+17,\\
        18k+7 &= 18(k-1)+6 + 19.
    \end{align*}
    For $k\geq 0$, we have
    \begin{align*}
        18k+8 &=  18k+6 + 2, \\
        18k+9 &= 18k+6 + 3,\\
        18k+11 &= 18k+ 6 + 5, \\
        18k+13 &=  18+ 6 + 7,\\
        18k+14 &= 18k+12 + 2,\\
        18k+15 &= 18k+12+3,\\
        18k+17 &= 18k+ 12 + 5.
    \end{align*}
    So, natural numbers $n$ congruent to $1$, $5$, $7$, $8$, $9$, $11$, $13$, $14$, $15$, or $17$ modulo $18$ can be written as a sum of a Zumkeller number and a prime if $n \geq 8$. Since the smallest Zumkeller number is $6$ and the smallest prime is $2$, $n$ cannot be written as a sum of a Zumkeller number and a prime if $n<6+2=8$.
\end{proof}

Finally, we prove that there are infinitely many natural numbers that cannot be expressed as a sum of a Zumkeller number and a square or as a sum of a Zumkeller number and a prime.

In the proceeding lemmas and theorems, we let $\sigma_{-1}(n)$ denote $\frac{\sigma(n)}{n}$. Note that for a prime $p$ and a natural number $a$,
\begin{equation}\label{eq: sigma_-1_inequality}
    \sigma_{-1}\left(p^a\right) = 1 + \frac{1}{p} + \cdots + \frac{1}{p^a} < \frac{p}{p-1}.
\end{equation}
Also, we denote primes in the ascending order as $p_1 =2$, $p_2=3$, $p_3=5$, and so on. Additionally, for every $\lambda > 1$ and natural number $k$, we define
$$u_{\lambda}(k) := \min \left\{s: \prod_{j=k+1}^s \frac{p_j}{p_j-1} \geq \lambda \right\},$$
and we define $A_{\lambda
}(k)$ to be the smallest natural number such that $\sigma_{-1}\left(A_{\lambda}(k)\right) \geq \lambda$ and $A_{\lambda}(k)$ is not divisible by $p_1, \cdots, p_k$. It is guaranteed that $u_{\lambda}(k)$ is well-defined since
$$\prod_{p \text{ prime, } p > p_k}\frac{p}{p-1}$$
diverges to infinity (see \cite[Section 22.7]{Hardy_Wright}).

Now we prove some lemmas which lead to the final results that we want to show. The proof of Lemma \ref{lem: from_abundant_paper} is similar to that of \cite[Lemma 1]{abundant_paper} with some small changes.

\begin{lemma}\label{lem: from_abundant_paper}
    $A_{\lambda}(k)$ is divisible by $p_{k+1}p_{k+2}\cdots p_{u_\lambda(k)}$.
\end{lemma}
\begin{proof}
    Let $M = p_{k+1}p_{k+2}\cdots p_{u_\lambda(k)}$. For the sake of contradiction, suppose $M\nmid A_{\lambda}(k)$. Let the prime factorization of $A_\lambda(k)$ be
    $$A_\lambda(k) = \prod_{i=1}^t q_i^{a_i}$$
    for distinct primes $q_1 < q_2 < \cdots < q_t$ and natural numbers $a_i$. Since $A_\lambda(k)$ is not divisible by $p_1, \cdots, p_k$, we have $q_1 > p_k$. Hence, $q_i\geq p_{k+i}$ for all $1\leq i \leq t$.
    
    We claim that $t \geq u_\lambda(k) - k$. For the sake of contradiction, assume $t<u_\lambda(k) - k$. Since $t\leq u_\lambda(k) - k - 1$ and $q_i\geq p_{k+i}$,
    \begin{align*}
        \sigma_{-1}(A_\lambda(k)) &= \sigma_{-1}\left(q_1^{a_1}\right)\cdots \sigma_{-1}\left(q_t^{a_t}\right) \\
        &\leq \sigma_{-1}\left(p_{k+1}^{a_1}\right)\cdots\sigma_{-1}\left(p_{k+t}^{a_{t}}\right)\\
        &\leq \sigma_{-1}\left(p_{k+1}^{a_1}\right)\cdots\sigma_{-1}\left(p_{u_\lambda(k) -1 }^{a_{u_\lambda(k)-k-1}}\right).
    \end{align*}
    Then, by equation (\ref{eq: sigma_-1_inequality}) and the definition of $u_{\lambda}(k)$, we have
    $$\sigma_{-1}(A_\lambda(k)) < \frac{p_{k+1}}{p_{k+1}-1}\cdots\frac{p_{u_{\lambda}(k)-1}}{p_{u_{\lambda}(k)-1} -1} < \lambda,$$
    which contradicts the fact that $\sigma_{-1}(A_\lambda(k)) \geq \lambda$. Hence, $t\geq u_\lambda(k)-k$.

    Since $M \nmid A_\lambda(k)$ by our assumption, $p_j \nmid A_\lambda(k)$ for some $k+1\leq j \leq u_\lambda(k)$. This, combined with the fact that $t\geq u_\lambda(k)-k$, implies that at least one of the primes $q_i$ dividing $A_{\lambda}(k)$ must be greater than $p_{u_\lambda(k)}$. Then, we have $q_t \geq q_i > p_{u_\lambda(k)} \geq p_j$. Since $q_t > p_j$ implies $\sigma_{-1}\left(p_j^{a_t}\right) > \sigma_{-1}\left(q_t^{a_t}\right)$, we have
    $$\sigma_{-1}\left(q_1^{a_1}\cdots q_{t-1}^{a_{t-1}} p_j^{a_t}\right) > \sigma_{-1}\left(q_1^{a_1}\cdots q_{t-1}^{a_{t-1}} q_t^{a_t}\right) \geq \lambda.$$
    Since
    $$q_1^{a_1}\cdots q_{t-1}^{a_{t-1}} p_j^{a_t} < q_1^{a_1}\cdots q_{t-1}^{a_{t-1}} q_t^{a_t} = A_\lambda(k)$$
    because $p_j < q_t$, this contradicts the minimality of $A_\lambda(k)$. Therefore, $A_\lambda(k)$ is divisible by $M$.
\end{proof}

\begin{lemma}\label{lem: sigma_at_least_2n}
    If $n$ is a Zumkeller number, then $\sigma(n) \geq 2n$.
\end{lemma}
\begin{proof}
    See \cite[Proposition 2(iii)]{Rao_Zumkeller_num} for proof.
\end{proof}

\begin{lemma}\label{lem: Zumkeller_lower_bound}
    If $s$ is a Zumkeller number such that $4\nmid s$ and $s$ is coprime to the first $k-1$ odd primes, then
    $$s \geq \min \left\{2 p_{k+1}\cdots p_{u_{4/3}(k)}, p_{k+1}\cdots p_{u_2(k)}\right\}.$$
\end{lemma}
\begin{proof}
    If $s$ is a Zumkeller number, then $\sigma_{-1}(s) \geq 2$ by Lemma \ref{lem: sigma_at_least_2n}. Since $4 \nmid s$ and $s$ is coprime to the first $k-1$ odd primes, $s=2^im$ for some non-negative integers $i$ and $m$ such that $i\in \{0,1\}$ and $m$ is not divisible by the first $k$ primes. Hence,
    $$\sigma_{-1}(s) = \sigma_{-1}\left(2^i\right)\sigma_{-1}(m) \geq 2,$$
    which implies that $\sigma_{-1}(m)\geq \frac{4}{3}$ if $i=1$ and that $\sigma_{-1}(m)\geq 2$ if $i=0$. By Lemma \ref{lem: from_abundant_paper},
    $$m \geq  p_{k+1}\cdots p_{u_{4/3}(k)}$$
    if $i=1$, and
    $$m \geq p_{k+1}\cdots p_{u_2(k)}$$
    if $i=0$. Therefore,
    $$s \geq \min \left\{2 p_{k+1}\cdots p_{u_{4/3}(k)}, p_{k+1}\cdots p_{u_2(k)}\right\}.$$
\end{proof}

\begin{lemma}\label{lem: k'_Zum_square}
    There exists an absolute constant $k'$ such that
    $$4p_2p_3\cdots p_k < \min\left\{2 p_{k+1}\cdots p_{u_{4/3}(k)}, p_{k+1}\cdots p_{u_2(k)}\right\}$$
    for all $k\geq k'$.
\end{lemma}
\begin{proof}
    By Mertens' third theorem,
    $$\prod_{p \leq x} \left(1-\frac{1}{p}\right)^{-1} \sim e^{\gamma} \ln{x},$$
    where $\gamma$ is the Euler-Mascheroni constant (see \cite[Theorem 429]{Hardy_Wright}). Then,
    $$\prod_{p_k<p\leq 3p_k}\left(1-\frac{1}{p}\right)^{-1} \sim \frac{\ln{3p_k}}{\ln{p_k}} \sim 1.$$
    Hence, for sufficiently large $k$, we have
    $$\prod_{p_k<p\leq 3p_k}\left(1-\frac{1}{p}\right)^{-1} < \frac{4}{3}.$$
    Thus, by the definition of $u_{4/3}(k)$ and $u_2(k)$, $p_{u_{4/3}(k)} > 3p_k$ and $p_{u_{2}(k)} > 3p_k$ for sufficiently large $k$.

    By the Prime Number Theorem (see \cite[Theorems 6 and 420]{Hardy_Wright}),
    $$\ln\left(4p_2\cdots p_k\right) \sim p_k,$$
    $$\ln\left(2p_{k+1}\cdots p_{u_{4/3}(k)}\right) \sim p_{u_{4/3}(k)} - p_k,$$
    and
    $$\ln\left(p_{k+1}\cdots p_{u_{2}(k)}\right) \sim p_{u_{2}(k)} - p_k.$$
    Since
    $$p_{u_{4/3}(k)} - p_k > 3p_k - p_k > p_k$$
    and
    $$p_{u_{2}(k)} - p_k > 3p_k - p_k > p_k$$
    for sufficiently large $k$, we have that
    $$\ln\left(4p_2\cdots p_k\right) < \ln\left(2p_{k+1}\cdots p_{u_{4/3}(k)}\right)$$
    and
    $$\ln\left(4p_2\cdots p_k\right) < \ln\left(p_{k+1}\cdots p_{u_{2}(k)}\right)$$
    for sufficiently large $k$. Therefore, there exists an absolute constant $k'$ such that for all $k\geq k'$,
    $$4p_2p_3\cdots p_k < \min\left\{2 p_{k+1}\cdots p_{u_{4/3}(k)}, p_{k+1}\cdots p_{u_2(k)}\right\}.$$
\end{proof}

\begin{theorem}\label{thrm: inf_no_Zum_square}
    There exist infinitely many natural numbers which cannot be written as a sum of a Zumkeller number and a square.
\end{theorem}
\begin{proof}
    For the sake of contradiction, assume that there are only finitely many natural numbers which cannot be written as a sum of a Zumkeller number and a square. Then, there exists a natural number $N$ such that $n$ can be written as a sum of a Zumkeller number and a square for all $n\geq N$.
    
    Since there are infinitely many primes of the form $4m+3$ (see \cite[Theorem 11]{Hardy_Wright}), we can choose a prime $p_j$ congruent to $3$ modulo $4$ such that $p_j \geq N+1$. Let $k$ be a natural number such that $k > \max\{k', j\}$, where $k'$ is defined in Lemma \ref{lem: k'_Zum_square}. By the Chinese Remainder Theorem, there exists a natural number $n\leq 4p_2\cdots p_k$ satisfying
    \begin{align*}
        n &\equiv r_1\, (\text{mod } 4),\\
        n &\equiv r_2\, (\text{mod } p_2),\\
        &\:\;\vdots\\
        n &\equiv -1\, (\text{mod } p_j),\\
        &\:\;\vdots\\
        n &\equiv r_k\, (\text{mod } p_k),
    \end{align*}
    where $r_1, r_2, \cdots, r_k$ are some quadratic non-residues modulo $4, p_2, p_3, \cdots, p_k$. Note that $r_1$ could be $2$ or $3$ and that $-1$ is a quadratic non-residue modulo $p_j$ (see \cite[Theorem 82]{Hardy_Wright}).

    Since $k > k'$ and $n\leq 4p_2\cdots p_k$, by Lemma \ref{lem: k'_Zum_square},
    $$n < \min\left\{2 p_{k+1}\cdots p_{u_{4/3}(k)}, p_{k+1}\cdots p_{u_2(k)}\right\}.$$
    Then, for any $0 < x < \sqrt{n}$,
    $$n-x^2 < \min\left\{2 p_{k+1}\cdots p_{u_{4/3}(k)}, p_{k+1}\cdots p_{u_2(k)}\right\}.$$
    Also, since $n$ satisfies the system of congruences described above, $n-x^2$ is not divisible by $4, p_2, p_3, \cdots, p_k$. Thus, by Lemma \ref{lem: Zumkeller_lower_bound}, $n-x^2$ is not a Zumkeller number for any $0<x<\sqrt{n}$. Hence, $n$ cannot be written as a sum of a Zumkeller number and a square.

    Since $n \equiv -1\,(\text{mod } p_j)$, $n+1\geq p_j \geq N+1$, which implies $n\geq N$. Thus, there exists a natural number $n\geq N$ such that $n$ cannot be written as a sum of a Zumkeller number and a square, which is a contradiction. Therefore, there exist infinitely many natural numbers which cannot be written as a sum of a Zumkeller number and a square.
\end{proof}

\begin{lemma}\label{lem: k'_Zum_prime}
    There exists an absolute constant $k''$ such that for all $k\geq k''$, if $n$ is coprime to the first $k$ primes and $\sigma_{-1}(n) \geq \frac{4}{3}$, then $n > p_1^2p_2^2\cdots p_k^2$.
\end{lemma}
\begin{proof}
    Assume that $n$ is coprime to the first $k$ primes and that $\sigma_{-1}(n) \geq \frac{4}{3}$. Then, by Lemma \ref{lem: from_abundant_paper},
    $$n \geq A_{4/3}(k) \geq p_{k+1}\cdots p_{u_{4/3}(k)}.$$
    By Mertens' Third Theorem,
    $$\prod_{p \leq x} \left(1-\frac{1}{p}\right)^{-1} \sim e^{\gamma} \ln{x},$$
    where $\gamma$ is the Euler-Mascheroni constant (see \cite[Theorem 429]{Hardy_Wright}). Then,
    $$\prod_{p_k<p\leq 4p_k}\left(1-\frac{1}{p}\right)^{-1} \sim \frac{\ln{4p_k}}{\ln{p_k}} \sim 1.$$
    Hence, for sufficiently large $k$, we have
    $$\prod_{p_k<p\leq 4p_k}\left(1-\frac{1}{p}\right)^{-1} < \frac{4}{3}.$$
    Thus, by the definition of $u_{4/3}(k)$, $p_{u_{4/3}(k)} > 4p_k$ for sufficiently large $k$. By the Prime Number Theorem (see \cite[Theorems 6 and 420]{Hardy_Wright}),
    $$\ln\left(p_1^2p_2^2\cdots p_k^2\right) \sim 2p_k,$$
    and
    $$\ln\left(p_{k+1}\cdots p_{u_{4/3}(k)}\right) \sim p_{u_{4/3}(k)} - p_k.$$
    Since
    $$p_{u_{4/3}(k)} - p_k > 4p_k - p_k > 2p_k$$
    for sufficiently large $k$, we have that
    $$\ln\left(p_{k+1}\cdots p_{u_{4/3}(k)}\right) > \ln\left(p_1^2p_2^2\cdots p_k^2\right)$$
    for sufficiently large $k$. Therefore, there exists an absolute constant $k''$ such that for all $k\geq k''$,
    $$n \geq p_{k+1}\cdots p_{u_{4/3}(k)} > p_1^2p_2^2\cdots p_k^2.$$
\end{proof}

\begin{theorem}\label{thrm: inf_no_Zum_prime}
    There exist infinitely many natural numbers which cannot be written as a sum of a Zumkeller number and a prime.
\end{theorem}
\begin{proof}
    Let $k$ be any natural number such that $k\geq k''$, where $k''$ is defined in Lemma \ref{lem: k'_Zum_prime}, and let $n = p_1^2p_2^2\cdots p_k^2$. Consider $n-p$ for any prime $p$ less than $n$.

    If $p > p_k$, then $n-p$ is coprime to the first $k$ primes. Since $n-p < n = p_1^2p_2^2\cdots p_k^2$, by Lemma \ref{lem: k'_Zum_prime}, $$\sigma_{-1}(n-p) < \frac{4}{3} < 2.$$
    Thus, by Lemma \ref{lem: sigma_at_least_2n}, $n-p$ is not a Zumkeller number.

    If $p\leq p_k$, then $n-p$ is divisible by $p$ but not by $p^2$. Hence, $\frac{n-p}{p}$ is coprime to the first $k$ primes, and thus
    $$\sigma_{-1}\left(\frac{n-p}{p}\right) < \frac{4}{3}$$
    by Lemma \ref{lem: k'_Zum_prime} because
    $$\frac{n-p}{p} < p_1^2p_2^2\cdots p_k^2.$$
    Hence,
    $$\sigma_{-1}(n-p) = \sigma_{-1}(p)\sigma_{-1}\left(\frac{n-p}{p}\right) < \frac{3}{2}\cdot\frac{4}{3} = 2$$
    because $p$ and $\frac{n-p}{p}$ are coprime, and thus $n-p$ is not a Zumkeller number by Lemma \ref{lem: sigma_at_least_2n}.

    Therefore, there exist infinitely many natural numbers $n$ of the form $n = p_1^2p_2^2\cdots p_k^2$ that are not expressible as a sum of a Zumkeller number and a prime.
\end{proof}
Note that the natural numbers $n$ considered in the proof above are even since $p_1=2$. Hence, Theorem \ref{thrm: inf_no_Zum_prime} particularly shows that there are infinitely many even natural numbers which cannot be expressed as a sum of a Zumkeller number and a prime, and we conjecture the following.
\begin{conjecture}\label{conj: odd_Zum_prime}
   All odd natural numbers greater than $8$ and not equal to $21$ are expressible as a sum of a Zumkeller number and a prime.
\end{conjecture}

Proposition \ref{prop: Zum_plus_prime_mod_18} shows that all odd natural numbers greater than $8$ other than those of the form $18k+3$ are expressible as a sum of a Zumkeller number and a prime. Thus, proving Conjecture \ref{conj: odd_Zum_prime} would require proving that all natural numbers of the form $18k+3$ for some integer $k\geq 2$ are expressible as a sum of a Zumkeller number and a prime. We have computationally verified that $18k+3$ can be written as a sum of a Zumkeller number and a prime for all $2\leq k \leq 1000$, and we suspect this to be true.

Regarding the natural numbers not expressible as a sum of a Zumkeller number and a square or as a sum of a Zumkeller number and a prime, we also have the following conjecture.

\begin{conjecture}\label{conj: lim_conjecture}
    Let $nzs(x)$ and $nzp(x)$ denote the number of natural numbers less than or equal to $x$ not expressible as a sum of a Zumkeller number and a square and the number of natural numbers less than or equal to $x$ not expressible as a sum of a Zumkeller number and a prime, respectively. Then,
    $$\lim_{x \to \infty} \frac{nzs(x)}{x} =\lim_{x \to \infty} \frac{nzp(x)}{x} = 0.$$
\end{conjecture}

\section{Future Research}\label{sec: fut_research}

In Section \ref{sec: Zum_arith}, we proved a result regarding Zumkeller numbers in arithmetic progressions. As an extension of this result, one might study quadratic, cubic, and biquadratic representations of Zumkeller numbers. Also, as a continuation of the result regarding sums of two Zumkeller numbers proven in Section \ref{sec: additive_rep_zum}, one could also explore numbers that are representable as a sum of three or four Zumkeller numbers. Additionally, analogically to the results regarding other additive representations, it might be of interest to study numbers that are expressible as a sum of a Zumkeller number and a higher power such as a cube, a biquadratic, etc. or as a sum of a Zumkeller number and a number of other types such as a deficient number, a perfect number, etc. Finally, it is also of the authors' interest to see Conjectures \ref{conj: odd_Zum_prime} and \ref{conj: lim_conjecture} solved.

\section{Disclosure Statement}
There is no potential conflict of interest that could influence this paper.

\section{Data Availability Statement}
There is no data associated with this article.

\bibliography{bibliography.bib}
\bibliographystyle{abbrv}

\end{document}